# An analysis of two modifications of the Petersburg game

ANDERS MARTIN-LÖF *

August 2007

### Abstract

Two modifications of the Petersburg game are considered:

1. Truncation, so that the player has a finite capital at his disposal.
2. A cost of borrowing capital, so that the player has to pay interest on the capital needed. In both cases limit theorems for the total net gain are derived, so that it is easy to judge if the game is favourable or not.

*Keywords:* Petersburg game, free lunch, arbitrage possibility, limited capital, cost of capital, profitability.
2000 Mathematics Subject Classifications: 60F05, 60E07.

## 1 Introduction

The well known Petersburg game is performed by tossing a fair coin until heads turns up for the first time. Let $T$ be this time. It has a geometric distribution, that is $P(T = k) = 2^{-k}$, $k = 1, 2, \ldots$.

Paul's gain in one game is $X = 2^T$, and this quantity has an infinite expectation, so that the ordinary law of large numbers does not provide a recipie for what is a fair price for playing several sucessive independent such games. Also this game allows the well known doubling strategy: As long as tails has come up double the stakes. Then the total amount spent is $1 + 2 + 4 + \ldots + 2^{T-1} = 2^T - 1$, so the net gain is $2^T - (2^T - 1) = 1$ with probability one, and we have a money machine providing an arbitrage possibility. This 'paradox' of course depends on the fact that Paul has an unbounded amount of money available for free to use as stakes. In this paper we want to analyse two natural modifications of this game:

1. Truncation. Paul has only a finite amount of capital, e.g. $2^c$, available for the stakes. This means that he has to stop if $T > c$, which event has probability $2^{-c}$. Otherwise he gains one as before and can continue to play. Then the expected gain is $1(1 - 2^{-c}) - (2^c - 1)2^{-c}$, so now the game is fair as has been pointed out e.g. in [1].

2. Introduction of interest on the capital. Paul can borrow money for the stakes without limit, but he has to pay interest at a constant rate on the loans needed. In this case the expeted present value of the net gain in a single game is zero as has also been pointed out in [1], so again the game is fair. Now Paul can continue to play forever, and we can consider the present value of all future costs and gains. We prove a limit theorem for these suitably normalised when the rate of interest $d$ goes to zero and hence the discount factor $r = 1/(1+d)$ goes to one. This involves a limit distribution very similar to that found in [2] for the total gain $S_N = X_1 + \ldots + X_N$ in a large number of independent games as $N = 2^n$ and $n \to \infty$, and in that case a simple asymptotic formula for the tail of the limit distribution can be obtained.

---
*Postal address: Department of Mathematics, Stockholm University, SE-106 91 Stockholm. E-mail: andersml@math.su.se

## 2 The truncated Petersburg game

The analysis of this game is quite straightforward. In each game the probability of a gain before stopping is $1 - 2^{-c}$, and the probability of stopping before a gain is $2^{-c}$, so the number of gains $M_c$ before stopping has a geometric distribution with $P(M_c \geq m) = (1-2^{-c})^m$. The total net gain before stopping is $V_c = M_c - (2^c - 1)$, and we get an exponential limit distribution as $c \to \infty$:

**Theorem 2.1** *When $c \to \infty$ $(M_c)2^{-c}$ converges in distribution to $U$ having an exponential distribution with $P(U \geq u) = e^{-u}, u \geq 0$, so that $(V_c)2^{-c}$ converges to $U - 1$.*

**Proof**: $P((M_c)2^{-c} \geq u) = (1 - 2^{-c})^{u2^c} \to e^{-u}$ as $c \to \infty$, so $P(U \geq u) = e^{-u}$, and for the total gain we see that $(V_c)2^{-c}$ converges to $U - 1$.

## 3  The Petersburg game with interest

Let us first consider the present value of the gain in a single game without the doubling strategy. The duration of the game $T$ has a geometric distribution as before: $P(T = k) = 2^{-k}, k = 1, 2, \ldots$, but now the present value of the gain is $(2r)^T$, which has a finite expectation equal to $\sum_{k=1}^{\infty} 2^{-k}(2r)^k = \sum_{k=1}^{\infty} r^k = r/(1-r)$. If an infinite number of games are played the gains occur at times $T_1, T_2, \ldots$ forming a renewal process with $T_i - T_{i-1}$ having the same geometric distribution as $T$. The present value of the total gain is now $V(r) = \sum_{i=1}^{\infty} r^{T_{i-1}}(2r)^{T_i - T_{i-1}}$. We now want to derive an asymptotic distribution for $V(r)$ when $r \to 1$. This can be done in the following way: As in [2] we scale time by a factor $N = 2^n \to \infty$ and put $r = e^{-a/N}$ with $1 \leq a < 2$. Then the renewal process has a deterministic limit: $(1/N)T_{Nt} \Rightarrow 2t$ since $E(T) = 2$. In the above reference the random walk generated by the successive gains: $S_k = \sum_{i=1}^{k} X_i$ with $X_i = 2^{T_i - T_{i-1}}$ is considered, and it is shown that the following limit theorem holds: $(1/N)S_{Nt} - nt \Rightarrow S(t)$, where $S(t)$ is a Lévy process which can be represented as follows: $S(t) = \sum_k Z_k(t)2^k$, where $Z_k(t)$ for $k$ positive are independent Poisson processes with mean $2^{-k}$, and for $k$ nonpositive are centered such processes. Its characteristic function thus has the following representation: $E(e^{izS(t)}) = e^{tl(z)}$ with

$$l(z) = \sum_{k=-\infty}^{0} 2^{-k}(e^{iz2^k} - 1 - iz2^k) + \sum_{k=1}^{\infty} 2^{-k}(e^{iz2^k} - 1),$$

so the Lévy measure $L(dx)$ has masses $2^{-k}$ at the positions $x = 2^k$. From this follows that the two-dimensional random walk $(T_k, S_k, k \geq 1)$ obeys the following limit theorem: $((1/N)T_{Nt}, (1/N)S_{Nt} - nt) \Rightarrow (2t, S(t))$ as $n \to \infty$. This can be used to find a limit distribution for $V(r)$: Consider

$$\begin{aligned}
2(1-r)V(r) - n &= 2(1-r)V(r) - 2an \int_{t=0}^{\infty} e^{-2at} dt \\
&\approx (2a/N) \sum_{i=1}^{\infty} e^{-aT_i/N}(X_i - n) \\
&\approx 2a \int_{t=0}^{\infty} e^{-2at} dS(t) \equiv 2aU,
\end{aligned}$$

because $S(t) \approx \sum_{i=1}^{Nt}(X_i - n)/N$. It is therefore interesting to try to estimate the right tail of the distribution of $U$. This can in fact be done in a way similar to that used for the study of the distribution of $S(t)$ in [2] because of the fact that $U$ has a Lévy representation quite similar to that of $S(t)$ :

**Lemma 3.1** $U$ has a Lévy distribution defined by $E(e^{izU}) = e^{g(z)/2a}$ with

$$g(z) = (2iz) + \sum_{l=-\infty}^{\infty} 2^{-l} \int_{2^l}^{2^{l+1}} (e^{izx} - 1 - izxc_l) dx/x,$$

where the centerings $c_l = 0$ for $l > 0$ and $c_l = 1$ for $l \leq 0$. The Lévy measure is hence defined by $L(dx) = 2^{-l} dx/2ax$ when $2^l < x < 2^{l+1}$ and $l$ is integer. This means that $U$ can be represented as $U = (2/2a) + \sum_{l=-\infty}^{\infty} W_l 2^l$, where $W_l$ are independent compound Poisson variables with

$$E(e^{izW_l}) = \exp\left((2^{-l}/2a) \int_1^2 (e^{izx} - 1 - izxc_l) dx/x\right).$$

**Proof:** Consider first the contribution from $Z_k(t)$ to $U = \int_0^\infty e^{-2at} dS(t)$. Since the increments $Z_k(dt)$ are independent its characteristic function is given by

$$\exp\left(\int_0^\infty (e^{iz2^k e^{-2at}} - 1 - iz2^k e^{-2at} c_k) 2^{-k} dt\right) =$$
$$\exp\left((2^{-k}/2a) \int_0^{2^k} (e^{izx} - 1 - izxc_k) dx/x\right),$$

with $x = 2^k e^{-2at}$. If we split the interval of integration into successive parts $(2^l, 2^{l+1})$ we get

$$\exp\left((2^{-k}/2a) \sum_{l=-\infty}^{k-1} \int_{2^l}^{2^{l+1}} (e^{izx} - 1 - izxc_k) dx/x\right).$$

Finally, summing also over $k$ we get

$$E(e^{izU}) = \exp\left(\sum_l \sum_{k=l+1}^{\infty} (2^{-k}/2a) \int_{2^l}^{2^{l+1}} (e^{izx} - 1 - izxc_k) dx/x\right).$$

For fixed $l$ we have $\sum_{k=l+1}^{\infty} 2^{-k} = 2^{-l}$ and $\sum_{k=l+1}^{\infty} 2^{-k} c_k = 0$ if $l \geq 0$, and $= \sum_{k=l+1}^{0} 2^{-k} c_k = (2^{-l} - 1) = 2^{-l}(1 - 2^l)$ if $l \leq 0$. Hence we see that

$$
\begin{aligned}
E(e^{izU}) &= \exp\left( (1/2a) \sum_{l>0} 2^{-l} \int_{2^l}^{2^{l+1}} (e^{izx} - 1) dx/x + \right. \\
&\qquad \left. (1/2a) \sum_{l \leq 0} 2^{-l} \int_{2^l}^{2^{l+1}} (e^{izx} - 1 - (izx)(1-2^l)) dx/x \right) \\
&= e^{g(z)/2a}
\end{aligned}
$$

with $g(z) =$

$$
(2iz) + \sum_{l>0} 2^{-l} \int_{2^l}^{2^{l+1}} (e^{izx} - 1) dx/x + \sum_{l \leq 0} 2^{-l} \int_{2^l}^{2^{l+1}} (e^{izx} - 1 - (izx)) dx/x.
$$

The decomposition $U = (2/2a) + \sum_l 2^l W_l$ is seen directly from this representation of $g(z)$ as well as the formula for the Lévy measure.

This Lévy representation is useful for getting an estimate of the right tail of the distribution of $U$ quite analogous to Theorem 3 in [2]. Let us outline the derivation of this. We first note that $g(z)$ is 'quasi-semi-stable' in the terminology of Lévy, i.e.

$$
\begin{aligned}
g(z2^{-m}) &= (2iz2^{-m}) + \sum_l 2^{-l} \int_{2^{l-m}}^{2^{l-m+1}} (e^{izx} - 1 - izxc_l) dx/x \\
&= (2iz2^{-m}) + 2^{-m} \sum_k 2^{-k} \int_{2^k}^{2^{k+1}} (e^{izx} - 1 - izxc_{k+m}) dx/x \\
&= 2^{-m}(g(z) + (iz) \sum_k (c_k - c_{l+m})) \\
&= 2^{-m}(g(z) + (izm)).
\end{aligned}
$$

This means that the characteristic function of $U_m \equiv 2^{-m}(U - m/2a)$ is given by

$$
\begin{aligned}
f_m(z) &\equiv E(e^{iz2^{-m}(U-m/2a)}) \\
&= e^{(g(z2^{-m}) - (izm)2^{-m})/2a} \\
&= e^{2^{-m} g(z)/2a} \\
&\approx 1 - 2^{-m} g(z)/2a
\end{aligned}
$$

as $m \to \infty$. Hence $2^m(f_m(z) - 1) \to g(z)/2a$. The term on the left is the Lévy exponent of the compound Poisson distribution whose Lévy measure $L_m(dx) = 2^m P(U_m \in dx)$, and the term on the right has Lévy measure $L(dx)$. Using the continuity theorem for Lévy exponents proved e.g. in [3] we can conclude that the tails of the distributions converge: $\int_x^\infty L_m(dy) \to \int_x^\infty L(dy) \equiv \bar{L}(x)$. $\bar{L}(x2^k)$ with $1 \leq x < 2$ can easily be calculated from Lemma 1:

$$\begin{aligned}
\bar{L}(x2^k) &= 2^{-k}\left(\int_{x2^k}^{2^{k+1}} dy/y + \sum_{l=k+1}^\infty 2^{-l}\int_{2^l}^{2^{l+1}} dy/y\right)/2a \\
&= 2^{-k}(\log 2 - \log x + \log 2)/2a \\
&= 2^{-k}(2\log 2 - \log x)/2a.
\end{aligned}$$

Taking $k = 0$ we see that the following asymptotic formula is valid for the distribution of $U$:

**Theorem 3.1** *As $m \to \infty$ and $1 \leq x < 2$ we have*
$2^m P(U > x2^m + m/2a) \to (2\log 2 - \log x)/2a.$

## 4 Conclusions

Remember that $V(r)$ is the present value of the total gain in an infinite sequence of games. Theorem 3.1 then allows us to find an initial premium which covers this value with an approximate given risk level. Remember also that we have $U = (V(r)/N) - (n/2a)$ with $r = e^{-a/N}$ and $N = 2^n \to \infty$. Then $P(V(r) > v) = P(V(r)/N > v/N) \approx P(U > v/N)$ if we neglect $n/2a$. This is $\approx 2^{-m}(2\log 2 - \log x)/2a$ if $v/N = x2^m$. This gives the simple estimate

$$\begin{aligned}
P(V(r) > v) &\approx (N/v)(2\log 2 - \log x)(x/2a) \\
&\approx x(2\log 2 - \log x)/2(1-r)v,
\end{aligned}$$

which can be used as a guide for choosing $v$. If the gambler pays the premium $v$ at the beginning his gain at time $Nt$ is approximately $e^{at}(V(r) - v)$, so he gains or looses at an exponential rate.

Let us finally analyse the total gains in an infinite sequence of Petersburg games with doubling in each game. In a single game the present value at the start of the gain minus the losses is

$$(2r)^T - r(1 + (2r) + \ldots + (2r)^{T-1}) = (2r)^T - r((2r)^T - 1)/(2r - 1)$$
$$= r/(2r - 1) - (2r)^T(1 - r)/(2r - 1).$$

The expected value of this is

$$r/(2r - 1) - ((1 - r)/(2r - 1)) \sum_{k=1}^{\infty} 2^{-k}(2r)^k = r/(2r - 1) - r/(2r - 1) = 0,$$

so this game is now fair. The present value of the total gain is now as before

$$\tilde{V}(r) = \sum_{i=1}^{\infty} r^{T_{i-1}}(r/(2r - 1) - ((1 - r)/(2r - 1))(2r)^{T_i - T_{i-1}}).$$

In the asymptotic approximation when $r = e^{-a/N}$ we see that $(2r - 1) \approx 1$ and we get $\tilde{V}(r) = \sum_{i=1}^{\infty} e^{-aT_{i-1}/N} - (1 - r)V(r)$. The first term is asymptotically $N \int_0^{\infty} e^{-2at} dt = N/2a$. and the second is $aU$. This means that $\tilde{V}(r) > 0$ with high probability, since $U$ remains finite as $N \to \infty$. We can get an approximate estimate of the ruin probability $R \equiv P(\tilde{V}(r) < 0)$, i.e. of $P(U > N/2a^2)$ using the previous formula if we take $x2^m = N/2a^2$. with $1 \leq x < 2$. We then get

$$R \approx ((2a^2x)/N)(2\log 2 - \log x) \approx (1 - r)(2ax)(2\log 2 - \log x).$$

As a numerical illustration consider a game which is played once a day and for which the interest rate is 4.46 % per year. Then $(1 - r) = a/N = 0.0446/365 = (1.22)10^{-4} = 2^{-13}$, and $x = a = 1$.

# References

[1] Aase, K. (2001): On the St. Petersburg Paradox, *Scand.Actuarial J.*, 2001;1: 69 - 78. MR1834973

[2] Martin-Löf, A. (1985): A Limit Theorem which clarifies the Petersburg Paradox, *J.Appl.Prob.* **22**, 634 - 643. MR0799286

[3] Feller, W. (1971): An Introduction to Probability Theory and its Applications, vol.2, 2nd ed., ch.XVII. *Wiley.*